\documentclass{birkmult}

\usepackage[draft]{hyperref}
\usepackage{txfonts}
\usepackage{amssymb}
\usepackage{tikz}\usetikzlibrary{positioning}
\usepackage
{subfigure}

\title{Nelson's Logical Diagrams}
\author[Andrew Aberdein]{Andrew Aberdein}
\address{School of Arts \& Communication\\ Florida Institute of Technology\\ Melbourne FL\\USA.}
\email{aberdein@fit.edu}
\date{\today}

\begin{document}

\begin{abstract}
As is now well known, the generalization of the square of opposition to a hexagon 
was discovered independently by three philosophers in the early 1950s:
Paul Jacoby, 
Augustin Sesmat, and 
Robert Blanch\'e
\cite{Jaspers16}.
Much less well known is the earlier use 
of similar diagrams by the German philosopher Leonard Nelson (1882--1927). 
This paper explores the significance of Nelson's work for 
the history of the JSB hexagon 
and for logical geometry more widely. 
\end{abstract}

\thanks{I delivered an earlier version of this paper at the 7th World Congress on the Square of Opposition, held at KU Leuven in September 2022. I am grateful to the audience in general and in particular to Ryan Christensen for supplying a copy of a diagram I had not seen. 
I am also grateful to the editors and referees of this volume for their encouraging comments.}

\maketitle


\section{Who was Leonard Nelson?%
\protect\footnote{This section is based on a similar discussion in \cite[p.~455]{Aberdein17d}.}}
A fascinating if neglected figure, Leonard Nelson (1882--1927) 
upheld a minority tradition within post-Kantian German philosophy: 
Nelson saw himself as a posthumous disciple of another neglected German philosopher, Jakob Friedrich Fries (1773--1843). Almost alone amongst Kant's nineteenth-century interpreters, Fries stressed the methodological aspects of Kant's work over its metaphysical implications. Fries's influence faded fast after his death, and his work was little known when Nelson encountered it as a student. Nelson made it his life's work to remedy this oversight.
Initially, this focus may have imperilled Nelson's own career: he secured a position at G\"ottingen only with the intervention of the mathematician David Hilbert, and over the objections of Edmund Husserl and other members of the philosophy faculty \cite[p.~122]{Reid70}.
Nelson's early death limited his direct influence---and frustratingly prevented what might have been a fruitful interaction with the pioneers of analytic philosophy. 
However, Nelson had a profound effect on his immediate circle amongst whom his posthumous influence was enduring.
Nelson was a prolific author across a range of philosophical topics, from metaphysics to the philosophy of law, but his chief focus was the methodology of philosophical reasoning itself. Although some of his major works have been translated into English \cite{Nelson49,Nelson56}, this project was abandoned before completion.
More recently, a series of lectures which had been omitted from Nelson's collected works has been published in German and English translation \cite{Nelson11,Nelson16}.
These lectures on ``Typische Denkfehler in der Philosophie'', or ``Typical Errors of Thinking in Philosophy'', were delivered at G\"ottingen in 1921.  
As we shall see, Nelson's contribution to logical geometry began with earlier work, but it is perhaps most fully developed in these lectures.

\section{Nelson's diagrams}
\begin{figure}[htbp]
\begin{center}
\scriptsize
\begin{tikzpicture}[thick]
\draw (-5.5,3.75) -- (-9.5,1.5) -- (-9.5,2.5) -- (-5.5,0) -- (-1.5,2.5) -- (-1.5,1.5) -- cycle;
\node at (-5.5,3.75) [above,align=center] {Jedes Urteil\\ ist entweder logisch\\ oder empirisch.};
\node at (-9.5,1.5) [below,align=center] {Die geometrischen Axiome\\ stammen aus der Logik.};
\node at (-9.5,2.5) [above,align=center] {Die geometrischen Axiome\\ stammen nicht\\ aus der Erfahrung.};
\node at (-5.5,0) [below,align=center] {Die geometrischen Axiome\\ stammen weder aus der Erfahrung\\ noch aus der Logik.};
\node at (-1.5,2.5) [above,align=center] {Die geometrischen Axiome\\ stammen nicht\\ aus der Logik.};
\node at (-1.5,1.5) [below,align=center] {Die geometrischen Axiome\\ stammen aus der Erfahrung.};
\end{tikzpicture}
\caption{Nelson's reconstruction of Kant \cite{Nelson11}}
\label{geometrie}
\end{center}
\end{figure}
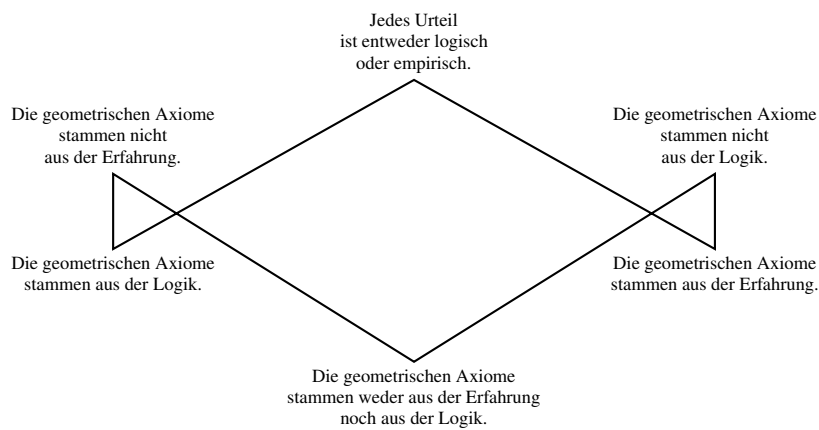

Nelson blamed many errors of philosophical reasoning on the fallacy of false dichotomy.
In particular, pairs of contrary propositions, which can be both false but not both true, are often treated as 
contradictory, so neither both true or both false.
Nelson used diagrams to depict this phenomenon, as in his reconstruction of Kant's defence of the synthetic a priori in Fig.~\ref{geometrie}.

\begin{figure}[htbp]
\begin{center}
\tiny
\begin{tikzpicture}[>=stealth,thick,every node/.style={draw,text width=3cm,align=flush center}] 
\node (PvQ) {Every judgment is either logical or empirical};
\node (notP) [below left=.5 of PvQ] {The axioms of geometry do not stem from experience};
\node (notQ) [below right=.5 of PvQ] {The axioms of geometry do not stem from logic};
\node (P) [below=1.5 of notQ] {The axioms of geometry stem from experience};
\node (Q) [below=1.5 of notP] {The axioms of geometry stem from logic};
\node (PnorQ) [below=3.5 of PvQ] {The axioms of geometry stem neither from experience nor from logic};
\draw[->]	(PvQ) -- (Q);
\draw[->]	(PvQ) -- (P);
\draw[->]	(notP) -- (Q);
\draw[->]	(notQ) -- (P);
\draw[->]	(notP) -- (PnorQ);
\draw[->]	(notQ) -- (PnorQ);
\end{tikzpicture}
\caption{Nelson's reconstruction of Kant, as presented by 
Leal 
\cite[p.~80]{Nelson16}}
\label{geometry}
\end{center}
\end{figure}
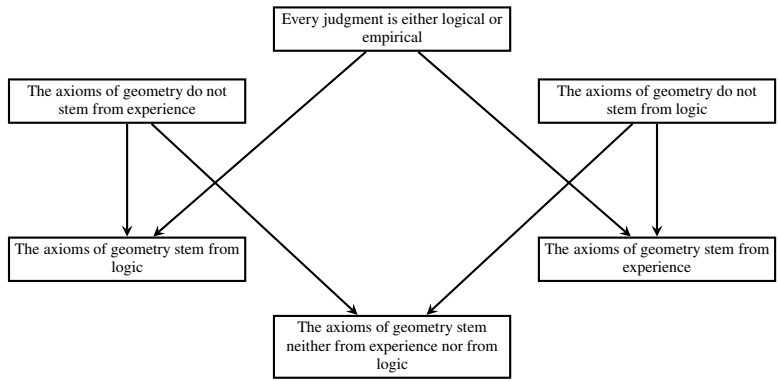

In the English edition of the 1921 G\"ottingen lectures, the translator and editor Fernando Leal adds arrowheads to clarify the direction of Nelson's inferences (Fig.~\ref{geometry}). Pairs of converging arrows are intended to be read as indicating arguments whose premisses are linked, thereby providing
joint support for the target proposition.
The starting positions of the two opposed parties are at upper left and right; at top is a shared presupposition that each relies on to derive conclusions at lower left and right, respectively. But Nelson (with Kant) urges rejection of this presupposition, leading to the bottom proposition, the conjunction of the ostensibly rival starting positions.

\begin{figure}[htbp]
\begin{center}
\subfigure[Leal's generalized Nelson diagram]{
\begin{tikzpicture}[>=stealth,thick,every node/.style={draw}] 
\node (PvQ) {$P\vee Q$};
\node (notP) [below left=.4 of PvQ] {$\neg P$};
\node (notQ) [below right=.4 of PvQ] {$\neg Q$};
\node (P) [below=of notQ] {$P$};
\node (Q) [below=of notP] {$Q$};
\node (PnorQ) [below=2.5 of PvQ] {$\neg(P\vee Q)$};
\draw[->]	(PvQ) -- (Q);
\draw[->]	(PvQ) -- (P);
\draw[->]	(notP) -- (Q);
\draw[->]	(notQ) -- (P);
\draw[->]	(notP) -- (PnorQ);
\draw[->]	(notQ) -- (PnorQ);
\end{tikzpicture}
\label{general}}
\quad\quad
\subfigure[A corresponding JSB hexagon]{
\begin{tikzpicture}[>=stealth,thick,every node/.style={draw}] 
\node (PvQ) {$P\vee Q$};
\node (P) [below right=.4 of PvQ] {$P$};
\node (Q) [below left=.4 of PvQ] {$Q$};
\node (notP) [below=.8 of Q] {$\neg P$};
\node (notQ) [below=.8 of P] {$\neg Q$};
\node (PnorQ) [below=2.5 of PvQ] {$\neg(P\vee Q)$};
\draw[<-]	(PvQ) -- (Q);
\draw[<-]	(PvQ) -- (P);
\draw[<-]	(notP) -- (Q);
\draw[<-]	(notQ) -- (P);
\draw[<-]	(notP) -- (PnorQ);
\draw[<-]	(notQ) -- (PnorQ);
\draw[-,red]	(PvQ) -- (PnorQ);
\draw[-,red]	(P) -- (notP);
\draw[-,red]	(Q) -- (notQ);
\draw[-,blue]	(P) -- (Q);
\draw[-,blue]	(P) -- (PnorQ);
\draw[-,blue]	(PnorQ) -- (Q);
\draw[-,green!75!black]	(PvQ) -- (notQ);
\draw[-,green!75!black]	(PvQ) -- (notP);
\draw[-,green!75!black]	(notP) -- (notQ);
\end{tikzpicture}
\label{hexagon}}
\caption{Leal's reconstruction of the general form of a Nelson diagram and its JSB 
counterpart}
\end{center}
\end{figure}
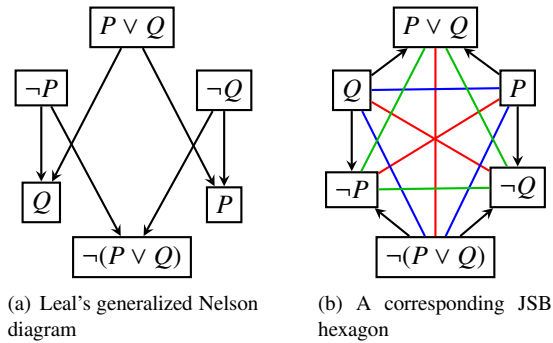

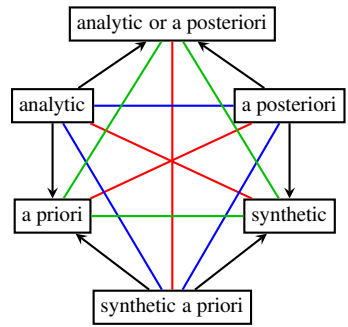
\begin{figure}[htbp]
\begin{center}
\footnotesize
\begin{tikzpicture}[>=stealth,thick,every node/.style={draw}] 
\node (PvQ) {analytic or a posteriori};
\node (Q) [below right=2.2em and -2em of PvQ.south east] {a posteriori};
\node (P) [below left=2.2em and -1.2em of PvQ.south west] {analytic};
\node (notP) [below=of Q] {synthetic};
\node (notQ) [below=of P] {a priori};
\node (PnorQ) [below=3.3 of PvQ] {synthetic a priori};
\draw[<-]	(PvQ) -- (Q);
\draw[<-]	(PvQ) -- (P);
\draw[<-]	(notP) -- (Q);
\draw[<-]	(notQ) -- (P);
\draw[<-]	(notP) -- (PnorQ);
\draw[<-]	(notQ) -- (PnorQ);
\draw[-,red]	(PvQ) -- (PnorQ);
\draw[-,red]	(P) -- (notP);
\draw[-,red]	(Q) -- (notQ);
\draw[-,blue]	(P) -- (Q);
\draw[-,blue]	(P) -- (PnorQ);
\draw[-,blue]	(PnorQ) -- (Q);
\draw[-,green!75!black]	(PvQ) -- (notQ);
\draw[-,green!75!black]	(PvQ) -- (notP);
\draw[-,green!75!black]	(notP) -- (notQ);
\end{tikzpicture}
\caption{B\'eziau's ``Kantian hexagon'' (after \cite[p.~27]{Beziau12})}
\label{Kantian}
\end{center}
\end{figure}

Although the diagrams in Nelson's 
1921 lectures are all example specific, 
Leal has proposed a reconstruction of their general form (Fig.~\ref{general}).
If $P$ and $Q$ are 
known to be contraries, then 
$P$, $Q$, and $\neg(P\vee Q)$ are pairwise contrary and 
$\neg P$, $\neg Q$, and $P\vee Q$ pairwise subcontrary. Of course each proposition is contradictory to its negation.
Hence, as I have observed elsewhere \cite[p.~458]{Aberdein17d}, the Nelson diagram in Fig.~\ref{general} is topologically equivalent to the JSB hexagon of opposition in Fig.~\ref{hexagon} 
(adding Jean-Yves B\'eziau's now standard colour-coding of contraries, subcontraries, and contradictories:
\cite[p.~3]{Beziau12}).
All that is required is to ``untwist'' the sides of the diagram and reverse the direction of Leal's arrows such that they now represent subalternation.
In particular, applying this procedure to Fig.~\ref{geometry} yields B\'eziau's ``Kantian Hexagon'' (Fig.~\ref{Kantian}).

\section{Further examples}
\begin{figure}[htbp]
\begin{center}
\footnotesize
\begin{tikzpicture}[>=stealth,thick,every node/.style={draw,text width=2.5cm,align=flush center}] 
\node (PvQ) {ISA acceptable $\vee$\\ WS necessary};
\node (Q) [below right=3em and 2em of PvQ.south east] {ISA acceptable};
\node (P) [below left=3em and 2em of PvQ.south west] {WS necessary};
\node (notP) [below=of Q] {WS not necessary};
\node (notQ) [below=of P] {ISA not acceptable};
\node (PnorQ) [below=3.5 of PvQ] {ISA not acceptable $\wedge$ WS not necessary};
\draw[<-]	(PvQ) -- (Q);
\draw[<-]	(PvQ) -- (P);
\draw[<-]	(notP) -- (Q);
\draw[<-]	(notQ) -- (P);
\draw[<-]	(notP) -- (PnorQ);
\draw[<-]	(notQ) -- (PnorQ);
\draw[-,red]	(PvQ) -- (PnorQ);
\draw[-,red]	(P) -- (notP);
\draw[-,red]	(Q) -- (notQ);
\draw[-,blue]	(P) -- (Q);
\draw[-,blue]	(P) -- (PnorQ);
\draw[-,blue]	(PnorQ) -- (Q);
\draw[-,green!75!black]	(PvQ) -- (notQ);
\draw[-,green!75!black]	(PvQ) -- (notP);
\draw[-,green!75!black]	(notP) -- (notQ);
\end{tikzpicture}
\caption{Nelson's political dilemma}
\label{political+}
\end{center}
\end{figure}
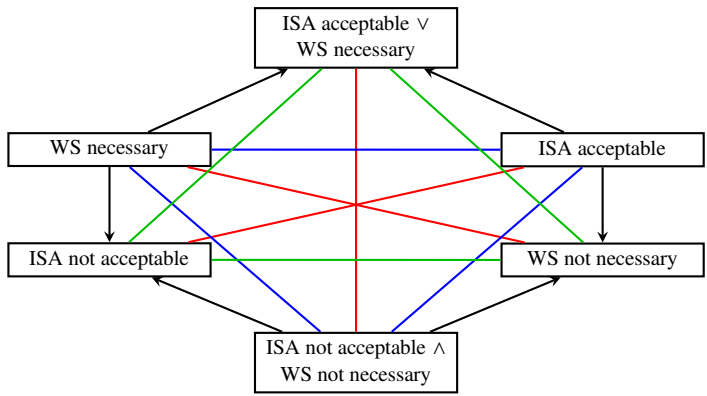

Nelson's 1921 lectures were not the first occasion on which he employed these distinctive logical diagrams. Other examples may be found in several of his earlier works (e.g.\ \cite[pp.~112, 123, 235, 372]{Nelson08}; \cite[pp.~22, 116, 282, 300, 308, 310, 322, 393]{Nelson17}).
For instance, Fig.~\ref{politische}, from 1915, analyses a political dilemma between inter-state anarchy and a world state (\cite[p.~22]{Nelson15}, reprinted as \cite[p.~22]{Nelson17})
|a question of more than conjectural interest when posed at the height of World War I. 
This example reoccurs in the 1921 lectures in somewhat simplified form as Fig.~\ref{political}.
I have recast it as a JSB hexagon in Fig.~\ref{political+}.
An even pithier political dilemma 
may be found in the English edition of one of Nelson's posthumously published works \cite[p.~252]{Nelson56}.

\begin{figure}[htbp]
\begin{center}
\tiny
\begin{tikzpicture}[thick]
\draw (-5.5,4.75) -- (-9.5,1.5) -- (-9.5,3.5) -- (-5.5,0) -- (-1.5,3.5) -- (-1.5,1.5) -- cycle;
\node at (-5.5,4.75) [above,align=center] {
Vollst\"andigkeit der Disjunktion von Staatenanarchie und Weltstaat.\\ (Jedes unabh\"angige politische Gemeinwesen\\ ist ein solches von Individuen.)};
\node at (-9.5,1.5) [below,align=center] {Behauptung der Notwendigkeit des Weltstaates.\\ (Es ist geboten, ein alle Individuen als\\ Glieder enthaltendes politisches\\ Gemeinwesen zu errichten.)};
\node at (-9.5,3.5) [above,align=center] {Behauptung der Verwerflichkeit der Staatenanarchie.\\ (Es ist geboten, da\ss\ sich die Staaten einem\\ sie umfassenden (Welt-)Gemeinwesen\\ unterwerfen.)};
\node at (-5.5,0) [below,align=center] {Unvollst\"andigkeit der Disjunktion von Staatenanarchie und Weltstaat.\\ (Ein unabh\"angiges politisches Gemeinwesen ist\\ nicht notwendig ein solches von Individuen.)};
\node at (-1.5,3.5) [above,align=center] {Bestreitung der Notwendigkeit des Weltstaates.\\ (Es ist nicht geboten, ein alle Individuen\\ als Glieder enthaltendes politisches\\ Gemeinwesen zu errichten.)};
\node at (-1.5,1.5) [below,align=center] {Bestreitung der Verwerflichkeit der Staatenanarchie.\\ (Es ist nicht geboten, da\ss\ sich die Staaten\\ einem sie umfassenden (Welt-)Gemeinwesen\\ unterwerfen.)};
\end{tikzpicture}
\caption{Nelson's political dilemma \cite[p.~22]{Nelson15}.}
\label{politische}
\end{center}
\end{figure}

\begin{figure}[htbp]
\begin{center}
\tiny
\begin{tikzpicture}[>=stealth,thick,every node/.style={draw,text width=3cm,align=flush center}] 
\node (PvQ) {Complete Disjunction between Inter-State Anarchy and the World-State};
\node (notP) [below left=.5 of PvQ] {Unacceptability of Inter-State Anarchy};
\node (notQ) [below right=.5 of PvQ] {Denial that a World-State is Necessary};
\node (P) [below=1.5 of notQ] {Denial that Inter-State Anarchy is Unacceptable};
\node (Q) [below=1.5 of notP] {Necessity of\\ a World-State};
\node (PnorQ) [below=3.5 of PvQ] {The above disjunction is incomplete};
\draw[->]	(PvQ) -- (Q);
\draw[->]	(PvQ) -- (P);
\draw[->]	(notP) -- (Q);
\draw[->]	(notQ) -- (P);
\draw[->]	(notP) -- (PnorQ);
\draw[->]	(notQ) -- (PnorQ);
\end{tikzpicture}
\caption{Nelson's political dilemma \cite[p.~173]{Nelson16}}
\label{political}
\end{center}
\end{figure}

\section{An imperfect example}
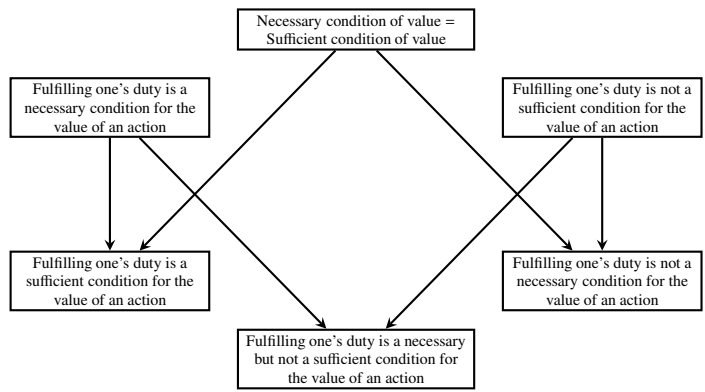
\begin{figure}[htbp]
\begin{center}
\tiny
\begin{tikzpicture}[>=stealth,thick,every node/.style={draw,text width=2.5cm,align=flush center}] 
\node (PvQ) [text width=3cm] {Necessary condition of value $=$ Sufficient condition of value};
\node (notP) [below left=.5 of PvQ] {Fulfilling one's duty is a necessary condition for the value of an action};
\node (notQ) [below right=.5 of PvQ] {Fulfilling one's duty is not a sufficient condition for the value of an action};
\node (P) [below=1.5 of notQ] {Fulfilling one's duty is not a necessary condition for the value of an action};
\node (Q) [below=1.5 of notP] {Fulfilling one's duty is a sufficient condition for the value of an action};
\node (PnorQ) [below=3.7 of PvQ,text width=3cm] {Fulfilling one's duty is a necessary but not a sufficient condition for the value of an action};
\draw[->]	(PvQ) -- (Q);
\draw[->]	(PvQ) -- (P);
\draw[->]	(notP) -- (Q);
\draw[->]	(notQ) -- (P);
\draw[->]	(notP) -- (PnorQ);
\draw[->]	(notQ) -- (PnorQ);
\end{tikzpicture}
\end{center}
\caption{Nelson on duty and value \cite[p.~186]{Nelson16}}
\label{duty}
\end{figure}

\begin{figure}[htbp]
\begin{center}
\footnotesize
\begin{tikzpicture}[>=stealth,thick,every node/.style={draw}] 
\node (PvQ) {$(V\rightarrow D)\rightarrow(D\rightarrow V)$};
\node (Q) [below right=2.2em and -1.6em of PvQ.south east] {$\neg(V\rightarrow D)$};
\node (P) [below left=2.2em and -1.2em of PvQ.south west] {$(D\rightarrow V)$};
\node (notP) [below=of Q] {$\neg(D\rightarrow V)$};
\node (notQ) [below=of P] {$(V\rightarrow D)$};
\node (PnorQ) [below=3.3 of PvQ] {$(V\rightarrow D)\wedge\neg(D\rightarrow V)$};
\draw[<-]	(PvQ) -- (Q);
\draw[<-
]	(PvQ) -- (P);
\draw[<-]	(notP) -- (Q);
\draw[<-]	(notQ) -- (P);
\draw[<-
]	(notP) -- (PnorQ);
\draw[<-]	(notQ) -- (PnorQ);
\draw[-,red]	(PvQ) -- (PnorQ);
\draw[-,red]	(P) -- (notP);
\draw[-,red]	(Q) -- (notQ);
\draw[-,blue]	(P) -- (Q);
\draw[-,blue]	(P) -- (PnorQ);
\draw[-,blue]	(PnorQ) -- (Q);
\draw[-,green!75!black]	(PvQ) -- (notQ);
\draw[-,green!75!black]	(PvQ) -- (notP);
\draw[-,green!75!black]	(notP) -- (notQ);
\end{tikzpicture}
\caption{Nelson on duty and value}
\label{duty+}
\end{center}
\end{figure}
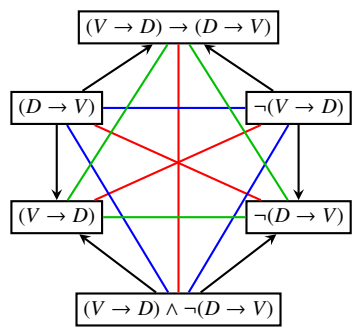

One limitation on this reconstruction is that Nelson does not always require that $P$ and $Q$ are contraries (despite the assumption made in Leal's generalization).
Consider, for example, Fig.~\ref{duty}, his analysis of the relationship of duty and value \cite[p.~186]{Nelson16}.
I have reconstructed it as Fig.~\ref{duty+}, where $D$ is ``one's duty is fulfilled'' and $V$ is ``one's action has value''. (The shared disjunctive premiss, here expressed as a conditional, has been simplified from Nelson's biconditional since its converse is not required for the inferences that he draws from it.)

However, we may observe two problems with Fig.~\ref{duty+}. Firstly, the shared disjunctive premiss, 
$(V\rightarrow D)\rightarrow(D\rightarrow V)$, is actually equivalent to one of the other statements, $(D\rightarrow V)$. Obviously, their contradictories, $(V\rightarrow D)\wedge\neg(D\rightarrow V)$ and $\neg(D\rightarrow V)$, are also equivalent. Hence the hexagon collapses into a square.
More importantly, the two ostensibly contrary statements $(D\rightarrow V)$ and $\neg(V\rightarrow D)$ are not contrary: they are both true when $D$ is false and $V$ is true, that is in a situation in which one's duty is not fulfilled but one's action does have value. 
Of course, this assumes a truth-functional account of the conditional, which is arguably anachronistic in the interpretation of Nelson.

\section{More than six vertices}
\begin{figure}[htbp]
\begin{center}
\tiny
\begin{tikzpicture}[>=stealth,thick,every node/.style={draw,text width=3cm,align=flush center}] 
\node (PvQ) {Complete disjunction between intuition and reflection as possible sources of knowledge};
\node (notP) [below left=.7 of PvQ] {Our intuition is not intellectual};
\node (notQ) [below right=.7 of PvQ] {Reflection is empty};
\node (P) [below=3 of notQ] {Metaphysics stems from intellectual intuition};
\node (Q) [below=3 of notP] {Metaphysics stems from reflection};
\node (PnorQ) [below=5 of PvQ] {Metaphysics stems from an immediate nonintuitive knowledge};
\node (PnandQ) [below=1.7 of PvQ] {We have no metaphysical knowledge};
\node (PandQ) [below=2.5 of PvQ] {We do have metaphysical knowledge};
\draw[->]	(PvQ) -- (Q);
\draw[->]	(PvQ) -- (P);
\draw[->]	(notP) -- (Q);
\draw[->]	(notQ) -- (P);
\draw[->]	(notP) -- (PnorQ);
\draw[->]	(notQ) -- (PnorQ);
\draw[->]	(PvQ) -- (PnandQ);
\draw[->]	(notP) -- (PnandQ);
\draw[->]	(notQ) -- (PnandQ);
\draw[->]	(PandQ) -- (P);
\draw[->]	(PandQ) -- (Q);
\draw[->]	(PandQ) -- (PnorQ);
\end{tikzpicture}
\caption{Nelson on immediate nonintuitive knowledge \cite[p.~199]{Nelson16}}
\label{ink}
\end{center}
\end{figure}
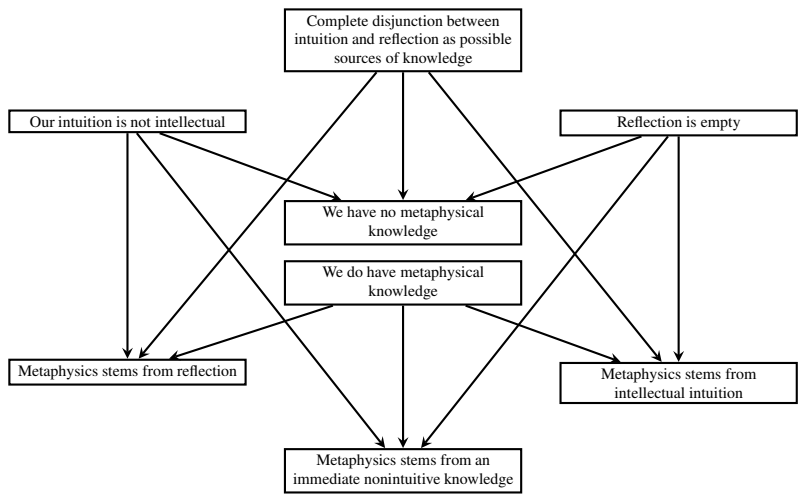

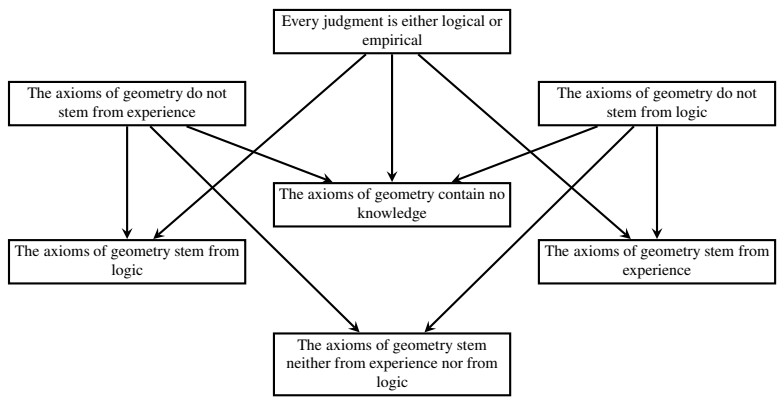
\begin{figure}[htbp]
\begin{center}
\tiny
\begin{tikzpicture}[>=stealth,thick,every node/.style={draw,text width=3cm,align=flush center}] 
\node (PvQ) {Every judgment is either logical or empirical};
\node (notP) [below left=.5 of PvQ] {The axioms of geometry do not stem from experience};
\node (notQ) [below right=.5 of PvQ] {The axioms of geometry do not stem from logic};
\node (P) [below=1.5 of notQ] {The axioms of geometry stem from experience};
\node (Q) [below=1.5 of notP] {The axioms of geometry stem from logic};
\node (PnorQ) [below=3.7 of PvQ] {The axioms of geometry stem neither from experience nor from logic};
\node (PnandQ) [below=1.7 of PvQ] {The axioms of geometry contain no knowledge};
\draw[->]	(PvQ) -- (Q);
\draw[->]	(PvQ) -- (P);
\draw[->]	(notP) -- (Q);
\draw[->]	(notQ) -- (P);
\draw[->]	(notP) -- (PnorQ);
\draw[->]	(notQ) -- (PnorQ);
\draw[->]	(PvQ) -- (PnandQ);
\draw[->]	(notP) -- (PnandQ);
\draw[->]	(notQ) -- (PnandQ);
\end{tikzpicture}
\caption{Nelson on conventionalism as ``Poincar\'e's false way out'' \cite[p.~97]{Nelson16}}
\label{Poincare}
\end{center}
\end{figure}

Perhaps the most dramatic departure that Nelson's diagrams make from the general form proposed by Leal is that some of them have more than six vertices.
Indeed, what appears to be the first such 
diagram in Nelson's work (Fig.~\ref{earliestD}) has eight vertices \cite[p.~57]{Nelson04}.
(The translation in Fig.~\ref{earliest} is taken from \cite[p.~146]{Nelson49}.)
It occurs at the end of a long footnote (spread over three pages in both German and English editions) in which Nelson proposes this diagram as an affirmative answer to a question posed by Kant: 
``Would it be possible to draw a priori a diagram of the history of philosophy, with which the eras of the various theories that philosophers have held (on the basis of the available data) would correspond as though they had had this very diagram before them and had proceeded in the knowledge of it?''%
\footnote{
This is Thomas K. Brown, III's translation from \cite[p.~144]{Nelson49}. 
Peter Heath, translating this passage for the Cambridge Edition of the Works of Immanuel Kant, renders the German ``ein Schema'' not as ``a diagram'', but less suggestively as ``a schema'' \cite[p.~418]{Kant02}.}
Nelson continues,
\begin{quotation}
Scientific interest in the history of philosophy is directed solely toward progress in the development of methods, not toward the results of individual philosophers---or at any rate toward these results only so far as they are dependent on the method followed. A law of the development of ideas can be formulated only in respect to method; and this law is rendered manifest by our diagram. This diagram is modeled on the organization of reason itself. The psychological point of view according to which it is constructed guarantees, on the one hand, that in all its subdivisions it is independent of standards that are either historically conditioned or arbitrarily assumed. Thus, it affords us a secure prescript by which we can peruse, in the light of principles, all the methodologically significant advances and errors in the history of the philosophical disciplines, in order to trace them back to their origin in reason itself.
(\cite[pp.~144 f., n. 11]{Nelson49}, translating \cite[pp.~55 f., n. 1]{Nelson04}).
\end{quotation}
A version of this diagram (Fig.~\ref{ink}) may 
be found in the 1921 lectures \cite[p.~199]{Nelson16}. 
The lectures also contain a seven-vertex diagram, Fig.~\ref{Poincare}, in which no additional premiss is required in order to reach the correct conclusion \cite[p.~97]{Nelson16}.

\begin{figure}[htbp]
\begin{center}
\tiny
\begin{tikzpicture}[every node/.style={text width=4cm,align=flush center},scale=.7]
  \coordinate (P1) at (.5,-4);
  \coordinate (P2) at (6.25,-.5);
  \coordinate (P3) at (12,-4);
  \coordinate (K1) at (12,-10.5);
  \coordinate (K2) at (6.25,-14);
  \coordinate (K3) at (.5,-10.5);
  \coordinate (Kk) at (6.25,-5);
  \coordinate (Dd) at (6.25,-9.5);
\node [above]  at (P1) {Faktische Pr\"amisse:\\ Wir besitzen Metaphysik.};
\node [text width=5cm,above] at (P2) {Dogmatische Pr\"amisse:\\ Alle Erkenntnis ist entweder Anschauung oder Reflexion.};
\node [above] at (P3) {Faktische Pr\"amisse:\\ Unsere Anschauung ist sinnlich.};
\node [text width=4.2cm,below]  at (K1) {Falsche Konsequens:\\ Also besitzen wir keine Metaphysik.\\ Konstitutives Princip: |\\Methodisches Princip: Unbegründbar. (Empirismus)};
\node [text width=7cm,below] at (K2) {Richtige Konsequenz:\\ Die Metaphysik entspringt aus nicht-anschaulicber unmittelbarer Erkenntnis.\\ Konstitutives Princip: Unmittelbare Erkenntnis der reinen Vernunft.\\ Methodisches Princip: Deduktion.\\ (Kriticismus)};
\node [text width=4.35cm,below] at (K3) {Falsche Konsequens:\\ Also besitzen wir intellektuelle Anschauung.\\ Konstitutives Princip: Intellektuelle Anschauung.\\ Methodisches Princip: Demonstration.\\ (Mysticismus)};
\node [below] at (Kk) {Falsche Konsequens:\\ Also entspringt die Metaphysik aus\\ der Reflexion.\\ Konstitutives Princip: Reflexion.\\ Methodisches Princip: Beweis.\\ (Logischer Dogmatismus)};
\node [above] at (Dd) {Faktische Pr\"amisse:\\Die reflektierte Erkenntnis ist mittelbar};
\draw [thick,rounded corners=1pt] (P2) -- (K3) -- (P1) -- (K2) -- (P3) -- (K1) -- cycle;
\draw [thick] (P1) -- (Kk) -- (P3);
\draw [thick] (K3) -- (Dd) -- (K1);
\draw [thick] (P2) -- (Kk);
\draw [thick] (K2) -- (Dd);
\end{tikzpicture}
\caption{Nelson's earliest published diagram \cite[p.~57]{Nelson04}.}
\label{earliestD}
\end{center}
\end{figure}
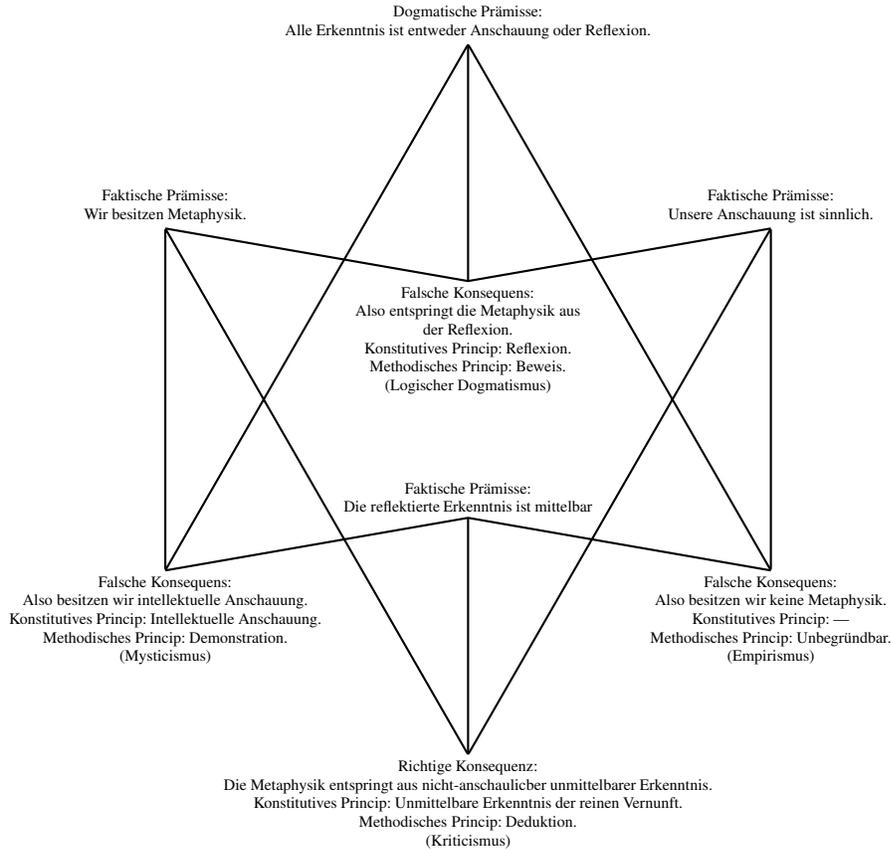

\begin{figure}[htbp]
\begin{center}
\scriptsize
\begin{tikzpicture}[every node/.style={text width=4cm,align=flush center},scale=.7]
  \coordinate (P1) at (.5,-4);
  \coordinate (P2) at (6.25,-.5);
  \coordinate (P3) at (12,-4);
  \coordinate (K1) at (12,-10.5);
  \coordinate (K2) at (6.25,-14);
  \coordinate (K3) at (.5,-10.5);
  \coordinate (Kk) at (6.25,-5);
  \coordinate (Dd) at (6.25,-9.5);
\node [above]  at (P1) {Factual premise: \\We possess metaphysics};
\node [text width=5cm,above] at (P2) {Dogmatic premise:\\ All knowledge is either intuition or reflection};
\node [above] at (P3) {Factual premise:\\ Our intuition is sensory};
\node [text width=4.2cm,below]  at (K1) {False consequence:\\ Therefore, we possess no metaphysics.\\ Constitutive principle: | \\Methodic principle: unverifiable. (Empiricism)};
\node [text width=7cm,below] at (K2) {Correct consequence:\\ Metaphysics has its source in nonintuitive immediate knowledge.\\ Constitutive principle: immediate knowledge of pure reason.\\ Methodic principle: deduction.\\ (Criticism)};
\node [text width=4.35cm,below] at (K3) {False consequence:\\ Therefore, we possess intellectual intuition.\\ Constitutive principle: intellectual intuition.\\ Methodic principle: demonstration. (Mysticism)};
\node [below] at (Kk) {False consequence:\\ Therefore, metaphysics has its source in reflection.\\ Constitutive principle: reflection.\\ Methodic principle: proof.\\ (Logical dogmatism)};
\node [above] at (Dd) {Factual premise:\\ Reflectional knowledge is derived};
\draw [thick,rounded corners=1pt] (P2) -- (K3) -- (P1) -- (K2) -- (P3) -- (K1) -- cycle;
\draw [thick] (P1) -- (Kk) -- (P3);
\draw [thick] (K3) -- (Dd) -- (K1);
\draw [thick] (P2) -- (Kk);
\draw [thick] (K2) -- (Dd);
\end{tikzpicture}
\caption{Nelson's earliest published diagram in translation \cite[p.~146]{Nelson49}.}
\label{earliest}
\end{center}
\end{figure}
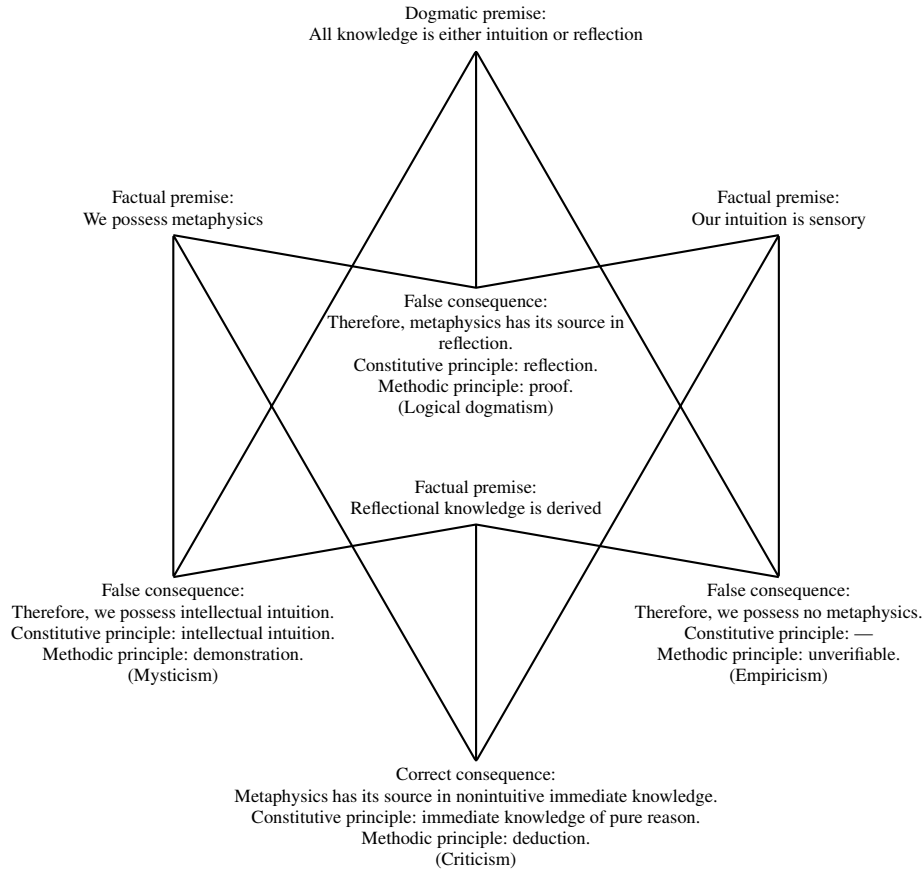

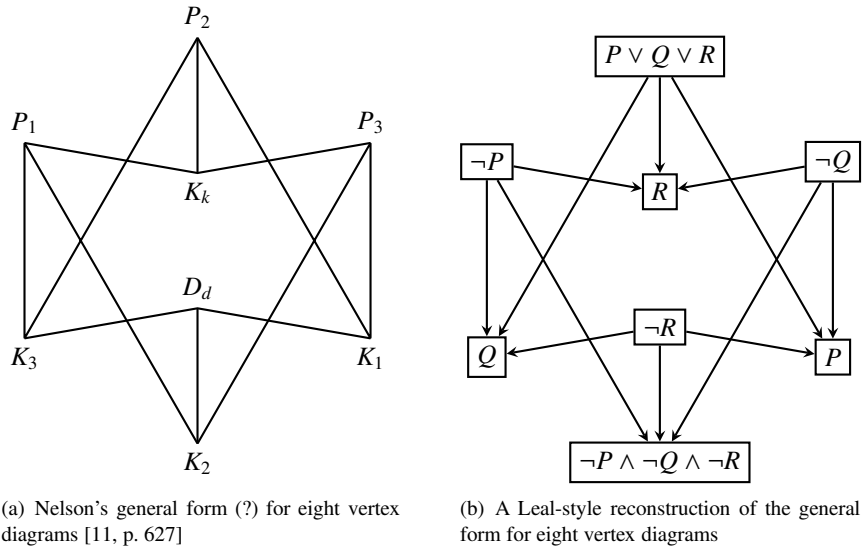
\begin{figure}[htbp]
\begin{center}
\subfigure[Nelson's general form (?) for eight vertex diagrams {\protect\cite[p.~627]{Nelson17}}]
{\begin{tikzpicture}[scale=.4]
  \coordinate (P1) at (.5,-4);
  \coordinate (P2) at (6.25,-.5);
  \coordinate (P3) at (12,-4);
  \coordinate (K1) at (12,-10.5);
  \coordinate (K2) at (6.25,-14);
  \coordinate (K3) at (.5,-10.5);
  \coordinate (Kk) at (6.25,-5);
  \coordinate (Dd) at (6.25,-9.5);
\node [above] at (P1) {$P_1$};
\node [above] at (P2) {$P_2$};
\node [above] at (P3) {$P_3$};
\node [below]  at (K1) {$K_1$};
\node [below] at (K2) {$K_2$};
\node [below] at (K3) {$K_3$};
\node [below] at (Kk) {$K_k$};
\node [above] at (Dd) {$D_d$};
\draw [thick,rounded corners=.5pt] (P2) -- (K3) -- (P1) -- (K2) -- (P3) -- (K1) -- cycle;
\draw [thick] (P1) -- (Kk) -- (P3);
\draw [thick] (K3) -- (Dd) -- (K1);
\draw [thick] (P2) -- (Kk);
\draw [thick] (K2) -- (Dd);
\end{tikzpicture}
\label{general8}}
\quad\quad
\subfigure[A Leal-style reconstruction of the general form for eight vertex diagrams]
{\begin{tikzpicture}[scale=.4,>=stealth,thick,every node/.style={draw}]
\node (P1) at (.5,-4) {$\neg P$};
\node (P2) at (6.25,-.5) {$P\vee Q\vee R$};
\node (P3) at (12,-4) {$\neg Q$};
\node (K1) at (12,-10.5) {$P$};
\node (K2) at (6.25,-14) {$\neg P\wedge\neg Q\wedge\neg R$};
\node (K3) at (.5,-10.5) {$Q$};
\node (Kk) at (6.25,-5) {$R$};
\node (Dd) at (6.25,-9.5) {$\neg R$};
\draw [->] (P1) -- (Kk);
\draw [->] (P2) -- (Kk);
\draw [->] (P3) -- (Kk);
\draw [->] (P1) -- (K3);
\draw [->] (P2) -- (K3);
\draw [->] (Dd) -- (K3);
\draw [->] (P2) -- (K1);
\draw [->] (P3) -- (K1);
\draw [->] (Dd) -- (K1);
\draw [->] (P3) -- (K2);
\draw [->] (P1) -- (K2);
\draw [<-] (K2) -- (Dd);
\end{tikzpicture}\label{general8+}}
\caption{Eight vertex diagrams}
\end{center}
\end{figure}

Leal does not present an account of the general form of the eight (or seven) vertex diagram analogous to that in Fig.~\ref{general}.
But, crucially, we may observe that the ``false way out'' can be understood as a third disjunct which must also be rejected in order to reach the correct conclusion. Hence, where the six vertex Nelson diagram is an analysis of a false dichotomy, the eight vertex diagram may be understood as an analysis of a false trichotomy.
The eight vertex diagram is also where Nelson comes closest to presenting a general form for his diagrams (Fig.~\ref{general8}).
However, it should be noted that the labels are all given specific definitions in the immediately preceding text \cite[p.~626 f.]{Nelson17}. So, although the labels are chosen in a manner suggestive of generality ($P$ for ``Pr\"amisse'', $K$ for ``Konsequenz'', and so forth), there is no direct indication that Nelson understood them in this manner, rather than as a passing typographical convenience.
On the other hand, the significance that Nelson attributes to this diagram (as indicated in the passage quoted above) and the heterogeneity of his examples suggest that he did understand it as capturing a widespread pattern of thought.
If we accept that the extra vertices may be treated as a third disjunction and its negation, we may generalize Leal's account of the six vertex diagram to the eight vertex case, as in Fig.~\ref{general8+}.

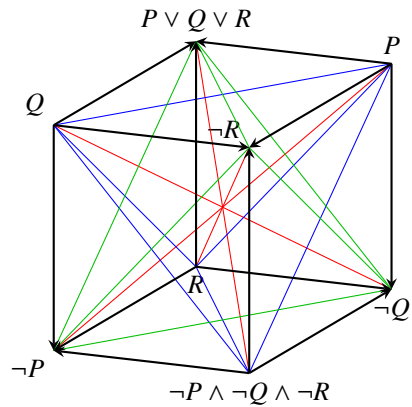
\begin{figure}[htbp]
\begin{center}
\begin{tikzpicture}[>=stealth,scale=3,rotate around y=-15]
  \coordinate (P1) at (0,0,1);
  \coordinate (P2) at (0,1,0);
  \coordinate (P3) at (1,0,0);
  \coordinate (K1) at (1,1,0);
  \coordinate (K2) at (1,0,1);
  \coordinate (K3) at (0,1,1);
  \coordinate (Kk) at (0,0,0);
  \coordinate (Dd) at (1,1,1);
\node [below left]  at (P1) {$\neg P$};
\node [above] at (P2) {$P\vee Q\vee R$};
\node [below] at (P3) {$\neg Q$};
\node [above]  at (K1) {$P$};
\node [below] at (K2) {$\neg P\wedge\neg Q\wedge\neg R$};
\node [above left] at (K3) {$Q$};
\node [below] at (Kk) {$R$};
\node [above left] at (Dd) {$\neg R$};
\draw[-,red]	(P2) -- (K2);
\draw[-,red]	(K1) -- (P1);
\draw[-,red]	(K3) -- (P3);
\draw[-,red]	(Kk) -- (Dd);
\draw[-,blue]	(K1) -- (K3);
\draw[-,blue]	(K1) -- (K2);
\draw[-,blue]	(K2) -- (K3);
\draw[-,blue]	(Kk) -- (K3);
\draw[-,blue]	(K1) -- (Kk);
\draw[-,blue]	(K2) -- (Kk);
\draw[-,green!75!black]	(P2) -- (P3);
\draw[-,green!75!black]	(P2) -- (P1);
\draw[-,green!75!black]	(P1) -- (P3);
\draw[-,green!75!black]	(Dd) -- (P3);
\draw[-,green!75!black]	(P2) -- (Dd);
\draw[-,green!75!black]	(P1) -- (Dd);
\draw[<-,thick]	(P2) -- (K3);
\draw[<-,thick]	(P2) -- (Kk);
\draw[<-,thick]	(P2) -- (K1);
\draw[<-,thick]	(P1) -- (K3);
\draw[<-,thick]	(P1) -- (Kk);
\draw[<-,thick]	(Dd) -- (K3);
\draw[<-,thick]	(P3) -- (K1);
\draw[<-,thick]	(P3) -- (Kk);
\draw[<-,thick]	(Dd) -- (K1);
\draw[<-,thick]	(P1) -- (K2);
\draw[<-,thick]	(Dd) -- (K2);
\draw[<-,thick]	(P3) -- (K2);
\end{tikzpicture}
\caption{A cube of opposition isomorphic to an eight vertex Nelson diagram}
\label{cube}
\end{center}
\end{figure}

I owe a final insight to a brief note by Kelley L. Ross, who observes that, since the eight vertex Nelson diagram also has twelve edges, it is isomorphic to a cube \cite{Ross06}.
Specifically, we may observe that it is isomorphic to the cube of opposition proposed by Alessio Moretti as resulting from four mutually exclusive propositions \cite[p.~25]{Moretti09a}.%
\footnote{Ross's cube differs from mine in some inessential details: I have oriented its vertices so as to make the isomorphism to Moretti's cube explicit.}
Thus we may reconstruct the eight vertex diagram as an instance of Moretti's cube of opposition (Fig.~\ref{cube}).
To justify this relationship, we generalize the assumption of contrariety made in mapping the six vertex diagram to the JSB hexagon to the claim that
$P$, $Q$, $R$, and $\neg P\wedge\neg Q\wedge\neg R$ are pairwise contrary. 
This in turn implies that 
$\neg P$, $\neg Q$, $\neg R$, and $P\vee Q\vee R$ are pairwise subcontrary. 
These comprise the blue and green tetrahedra, respectively, in Fig.~\ref{cube}.
All twelve edges represent implications: those from $\neg P\wedge\neg Q\wedge\neg R$ to one of its conjuncts or to $P\vee Q\vee R$ from one of its disjuncts are logically necessary; the other six follow from the assumption of contrariety.

\section{Conclusion}
It is the contention of this paper that the foregoing establishes Nelson as having anticipated what was once called the Blanch\'e hexagon and is now more often called the JSB hexagon. (Perhaps it should henceforth be called the NJSB hexagon.) Even more remarkably, Nelson can also be seen as anticipating Moretti's analogous cube of opposition.
In making these claims, we must acknowledge several important caveats.
Firstly, it is not entirely clear that Nelson was proposing a general scheme (tantalizing as Fig.~\ref{general8} may be). Even so, it is straightforward to generalize Nelson's work.
Secondly, Nelson's hexagons are not the convex hexagons familiar from Blanch\'e's work. But convex hexagons are not the only hexagons: Nelson's hexagons are 
self-intersecting hexagons.%
\footnote{Specifically, they are one of six types of self-intersecting hexagons with regular vertices, identified as the ``center-flip'' hexagon here: \url{https://en.wikipedia.org/wiki/Hexagon\#Self-crossing_hexagons}.}
It should also be noted that Jacoby's ``double triangle'' is not a convex hexagon either \cite[p.~44]{Jacoby50}.
Thirdly, some of Nelson's diagrams---indeed, the earliest published---have more than six vertices. But, as we have seen, this may be taken as suggesting Nelson was even further advanced in his anticipation of later work in logical geometry.%
\footnote{Although Lewis Carroll, for one, may be argued to have gone further and earlier \cite{Moretti14}.}
Lastly, but perhaps most importantly, there is no indication that Nelson understood his work as an extension of the traditional square of opposition, or in any way a contribution to the logic of categorical propositions. 
His diagrams make the entailment relations of the hexagon (or cube) explicit but do not represent the oppositional relations of contrariety, subcontrariety, or contradiction in graphical terms. Rather, these relations are set out in (or may be inferred from) the associated text. 
Relatedly, Nelson does not always seem to require that $P$ and $Q$ are contraries (despite the assumption made in Leal's generalization). Nonetheless, we might regard such exceptional cases as slips, rather than as evidence that Nelson was engaged in a different, incompatible project. 
If any or all of these concerns are granted, then Nelson's work in logical geometry should be seen only as an incomplete anticipation of that of Jacoby, Sesmat, and Blanch\'e. 
Even so, it would still warrant wider attention than it has so far received.

An even more speculative question is whether Nelson might have influenced Jacoby, Sesmat, or Blanch\'e. Perhaps the clearest presentation of Nelson's diagrams is in the 1921 lectures, which remained unpublished until long after the deaths of all three men. (Similarly, Lewis Carroll's work on logical geometry existed only as unpublished fragments until long after its contents had been independently rediscovered: \cite[p.~387]{Moretti14}.)
Although Nelson used his diagrams in many 
works published 
in German from early in the twentieth century, 
these publications do not seem to have had a wide reception even in Germany, so it seems quite likely that neither Jacoby, Sesmat, nor Blanch\'e would have encountered them. But, more significantly, in 1949 a major American publisher brought out an English translation of a selection of Nelson's work \cite{Nelson49}, 
including 
perhaps the earliest 
of Nelson's diagrams, shortly before Jacoby published his own logical diagrams \cite{Jacoby50}. This book was widely reviewed and it is certainly possible that Jacoby could have come across it. Even so, this is probably no more than a coincidence: Jacoby's presentation is quite different from Nelson's, and is directly concerned with traditional logic, as Nelson is not. The most likely conclusion would seem to be that 
Nelson's work 
in logical geometry was unknown to any of the later scholars he anticipated.


\begin{thebibliography}{10}
\providecommand{\url}[1]{{#1}}
\providecommand{\urlprefix}{URL }
\expandafter\ifx\csname urlstyle\endcsname\relax
  \providecommand{\doi}[1]{DOI~\discretionary{}{}{}#1}\else
  \providecommand{\doi}{DOI~\discretionary{}{}{}\begingroup
  \urlstyle{rm}\Url}\fi

\bibitem{Aberdein17d}
Aberdein, A.: Leonard {N}elson: A theory of philosophical fallacies.
\newblock Argumentation \textbf{31}(2), 455--461 (2017)

\bibitem{Beziau12}
B\'eziau, J.Y.: The power of the hexagon.
\newblock Log.\ Univers.\ \textbf{6}(1), 1--43 (2012)

\bibitem{Jacoby50}
Jacoby, P.: A triangle of opposites for types of propositions in {A}ristotelian
  logic.
\newblock New Scholasticism \textbf{24}(1), 32--56 (1950)

\bibitem{Jaspers16}
Jaspers, D., Seuren, P.A.: The square of opposition in {C}atholic hands: A
  chapter in the history of 20th-century logic.
\newblock Logique et Anal.\ \textbf{59}, 1--35 (2016)

\bibitem{Kant02}
Kant, I.: Theoretical Philosophy after 1781.
\newblock Cambridge University Press, Cambridge (2002).
\newblock Trans. Hatfield, G., Friedman, M.; Ed. Allison, H., Heath  P.

\bibitem{Moretti09a}
Moretti, A.: The geometry of standard deontic logic.
\newblock Log.\ Univers.\ \textbf{3}(1), 19--57 (2009)

\bibitem{Moretti14}
Moretti, A.: Was {L}ewis {C}arroll an amazing oppositional geometer?
\newblock Hist.\ Philos.\ Logic \textbf{35}(4), 383--409 (2014)

\bibitem{Nelson04}
Nelson, L.: Die kritische Methode und das Verh{\"a}ltnis der Psychologie zur
  Philosophie: Ein Kapitel aus der Methodenlehre.
\newblock Vandenhoeck \& Ruprecht, G\"ottingen (1904)

\bibitem{Nelson08}
Nelson, L.: {\"U}ber das sogenannte Erkenntnisproblem.
\newblock Vandenhoeck \& Ruprecht, G\"ottingen (1908)

\bibitem{Nelson15}
Nelson, L.: Ethische Methodenlehre.
\newblock Veit, Leipzig (1915)

\bibitem{Nelson17}
Nelson, L.: Kritik der praktischen Vernunft.
\newblock Veit, Leipzig (1917)

\bibitem{Nelson49}
Nelson, L.: Socratic Method and Critical Philosophy: Selected Essays.
\newblock Yale University Press, New Haven, CT (1949).
\newblock Trans. Brown, III, T.K.

\bibitem{Nelson56}
Nelson, L.: System of Ethics.
\newblock Yale University Press, New Haven, CT (1956).
\newblock Trans. Guterman, N.

\bibitem{Nelson11}
Nelson, L.: Typische Denkfehler in der Philosophie.
\newblock Felix Meiner, Hamburg (2011).
\newblock Ed. Brandt, A., Schroth, J.

\bibitem{Nelson16}
Nelson, L.: A Theory of Philosophical Fallacies.
\newblock Springer, Cham (2016).
\newblock Trans. Leal, F., Carus, D.

\bibitem{Reid70}
Reid, C.: Hilbert--Courant.
\newblock Springer, New York, NY (1970)

\bibitem{Ross06}
Ross, K.L.: On {N}elson's axiomatic diagrams.
\newblock The Proceedings of the Friesian School, Fourth Series  (2006).
\newblock \url{https://friesian.com/universl.htm#note-3}


\end{thebibliography}

\end{document}